# A Fast Solution Method for Large-scale Unit Commitment Based on Lagrangian Relaxation and Dynamic Programming

Jiangwei Hou, *Student Member, IEEE*, Qiaozhu Zhai, *Member, IEEE*, Yuzhou Zhou, *Member, IEEE*, Xiaohong Guan, *Fellow, IEEE*

*Abstract*—The unit commitment problem (UC) is crucial for the operation and market mechanism of power systems. With the development of modern electricity, the scale of power systems is expanding, and solving the UC problem is also becoming more and more difficult. To this end, this paper proposes a new fast solution method based on Lagrangian relaxation and dynamic programming. Firstly, the UC solution is estimated to be an initial trial UC solution by a fast method based on Lagrangian relaxation. This initial trial UC solution fully considers the system-wide constraints. Secondly, a dynamic programming module is introduced to adjust the trial UC solution to make it satisfy the unit-wise constraints. Thirdly, a method for constructing a feasible UC solution is proposed based on the adjusted trial UC solution. Specifically, a feasibility-testing model and an updating strategy for the trial UC solution are established in this part. Numerical tests are implemented on IEEE 24-bus, IEEE 118-bus, Polish 2383-bus, and French 6468-bus systems, which verify the effectiveness and efficiency of the proposed method.

*Index Terms*-- Fast solution, large-scale unit commitment, Lagrangian relaxation, dynamic programming.

## NOMENCLATURE

### A. Indices and sets

| | |
|---|---|
| $i$ | Index of units, $i = 1,2,...,N$ |
| $k$ | Index of iterations, $k = 1,2,...,K$ |
| $l$ | Index of transmission lines, $l = 1,2,...,L$ |
| $m$ | Index of load buses, $m = 1,2,...,M$ |
| $t$ | Index of scheduling periods, $t = 1,2,...,T$ |

### B. Variables

| | |
|---|---|
| $p$ | Power outputs of the units during the periods |
| $v_{l,t}^+, v_{l,t}^-, v_{0,t}^+, v_{0,t}^-$ | Non-negative scalar slack variables for the system-wide constraints. |
| $z$ | On/off status of the units during the periods. |
| $\hat{z}^*, \bar{z}^*$ | The trial UC solution and the new UC solution. |
| $\lambda, \hat{\lambda}^*, \hat{\lambda}_k^*$ | Dual variables, a trial dual solution, and an updated trial dual solution. |

### C. Parameters

| | |
|---|---|
| $a,b$ | Coefficients of the affine fuel cost functions. |
| $c_k$ | Step-size for updating $\hat{\lambda}_k^*$. |
| $d_{i,t}$ | Net load demand at bus $i$ and of period $t$ (MW). |
| $F_l$ | Transmission limit of transmission line $l$ (MW) |
| $\underline{P}, \overline{P}$ | Minimum/maximum generation capacities (MW). |
| $\beta_i^0$ | A constant threshold for determining $\hat{z}_i^*$. |
| $\Gamma$ | Power transfer distribution factors (PTDF). |
| $\Delta_i^+, \Delta_i^-$ | Ramp up/down limit of unit $i$ (MW). |
| $\varepsilon$ | Convergence tolerance of the method. |
| $\overline{\tau}_i, \underline{\tau}_i$ | Minimum up/down time of the unit $i$ (hour). |
| $\delta$ | Thresholds for optimizing the initial step size $c_0$. |
| $u$ | Bound of the initial step size $c_0$. |

### D. Functions

| | |
|---|---|
| $q(\cdot), \hat{q}(\cdot)$ | Lagrangian dual function of the UC problem, and an approximate Lagrangian dual function. |
| $g_{\hat{\lambda}_k^*}$ | The sub-derivative of $q(\lambda)$ at $\hat{\lambda}_k^*$. |
| $S(\cdot), C(\cdot)$ | Functions of the start-up costs and the fuel costs of units ($). |
| $\beta(\cdot)$ | Function for determining $\hat{z}^*$ with $\beta^0$. |

### E. Abbreviations

| | |
|---|---|
| trial UC | A trial UC solution. |
| new UC | A new UC solution that is adjusted from trial UC. |
| DPLR | The proposed method based on dynamic programming and Lagrangian relaxation. |
| NSTD | A new state transition diagram for the unit sub-problems of LR methods. |

## I. INTRODUCTION

UNIT commitment (UC) is one of the most critical problems focused on by the operators of power systems. The main objective of UC problems is to find a day-ahead scheduling plan for the on/off status of generating units. This scheduling plan achieves the minimum generation cost, and it guarantees the satisfaction of various system-wide and unit-wise constraints about generation, transmission, and system stability [1]-[3].

With the rapid expansion of power systems, the solving of UC problems becomes more and more difficult, and obtaining fast solutions for large-scale UC problems has become a pressing issue [4]. The main difficulty originates from that the mixed variables (binary and continuous) in the problem are heavily coupled by various constraints. In particular, the sys-

This work is supported by the Science and Technology Project of State Grid Corporation of China (5400-202199524A-0-5-ZN).
J. Hou, Q. Zhai, Y. Zhou, and X. Guan are with Systems Engineering Institute, MOEKLINNS Lab, Xi'an Jiaotong University, Xi'an 710049, China. (Corresponding author: qzzhai@sei.xjtu.edu.cn).

tem-wide constraints couple the variables of all the generating units at each scheduling period, and the unit-wise constraints couple the variables of a unit in adjacent periods [4], [5]. Indeed, this problem is proven to be NP-hard [4], and for decades it has been drawing much attention from academies and the electric power industry.

There have been various approaches for solving the UC problem, mainly including the mix-integer linear programming (MILP) [6]-[14], Lagrangian relaxation (LR) methods [15]-[19], dynamic programming [20], data-driven methods [21], [22], and heuristic methods [23]. These methods keep a proper balance between the efficiency and quality of the solutions, for UC problems of different sizes and characteristics. The heuristic methods, such as the genetic algorithm [23], are generally efficient in searching for locally optimum solutions, but their computation time is unpredictable; the data-driven methods improve the accuracy and efficiency based on prior knowledge of historical data [22]. To this end, this paper aims at providing a fast solution of high quality for large-scale UC problems without resorting to prior knowledge. The existing methods of this kind are mainly relevant to MILP methods [6]-[14] and LR methods [15]-[19].

The MILP methods enumeratively search for the optimal UC solution in the feasible domain of the binary variables with the branch-and-bound method, and they have nowadays been the researchers' preferred solution method mainly due to the improved performance of MILP solvers. The improvements along with their variants make it easier for the researchers to formulate and to solve MILP models with more modeling details. These details are constructed mainly to further improve the efficiency and/or the accuracy of MILP methods. For example, the tight MILP formulations in [6]-[9] are said to be more computationally efficient since their strengthened linear programs provide improved lower bounds of the objective; some other MILP formulations finer reflect the operational constraints of the systems that are usually coarsely approximated [10]. In addition, the Benders decomposition method and its accelerated versions for MILP also considerably improve the efficiency [11]-[14]. However, their computational complexity grows exponentially with the problem size, and it is reported that at times existing solvers cannot effectively solve the MILP problems within the time window for clearing the electricity market [15].

LR methods usually relax the coupling constraints and penalize violations of these constraints in the objective with Lagrange multipliers. This makes the problem separate into subproblems that can be efficiently solved. For example, the system-wide constraints are relaxed in [24], and this yields unit sub-problems for every unit. This sub-problem can be solved in polynomial time by some efficient dynamic programming methods [16], [17]. When the unit-wise constraints are further relaxed, the unit sub-problem can be further separated into unit-period sub-problems of different periods whose UC solution can be analytically calculated [15]. These UC solutions are found to be close to the optimal UC solution [15],[16], however, they are generally infeasible due to the relaxations. Thus, the solution and the associated multipliers are usually iteratively adjusted with sub-gradient methods to generate trial solutions converging to the optimal solution [16]. This approach usually has low computational complexity in practice, but it incurs the problems of zigzagging trial solutions since the sub-gradients are not always effective in updating the trial solutions [25]. And improvement of the convergence has perennially been an active research area. For example, the surrogate sub-gradient developed and utilized in [24], [26] improves the solution speed by not fully optimizing the relaxed problem; and in [27] the Lagrangian function is augmented with penalties of squared constraint violations to make the trial solutions converge faster.

It is observed from the above methods that the binary variables are the main cause of the computational difficulties. The MILP methods search for the optimal binary solution on a tree structure with branch-and-bound/cut methods, and thus have combinatorially explosive complexity. LR methods converge slowly, and the optimal LR solution (or even a feasible LR solution) can be computationally difficult to obtain for large-scale problems. The computational difficulty is partially relieved by the surrogate sub-gradient methods, but it is still unacceptable for large-scale problems [26]. Thus, it would be better to avoid searching for the solution from the tree structure and to guarantee the efficiency of sub-gradient methods in obtaining the optimal/feasible solution.

There do exist methods providing important insights for obtaining a fast feasible LR solution. Based on LR, an analytical function of Lagrange multipliers in [15] gives a fast trial UC solution that fully considers the satisfaction of the system-wide constraints but ignores the unit-wise constraints. In [17], a new state transition diagram (NSTD) of dynamic programming designed for the unit sub-problem, which only has hundreds of transition states, produces a fast UC solution strictly satisfying the unit-wise constraints. In addition, the feasibility of a trial UC solution can be efficiently examined with a linear program by the method in [28].

In this paper, the advantages of the just mentioned methods are combined to give a systematic method for providing a fast UC solution. Specifically, the fast trial solution of [15] is iteratively generated and efficiently adjusted with the NSTD in [17]. In this way, the method can generate a fast solution that is feasible to all the constraints. The adjustment to the iteratively generated trial solution is based on a constructed feasibility-testing problem. This problem efficiently yields information for improving the feasibility of the generated trial solutions. In addition, the proposed method searches for the binary solution in the NSTD with dynamic programming, and thus it has an extremely lower computational complexity than MILP. In numerical tests, the method finds the solution after a few trial solutions have been explored, and the problems of the 6468-bus system can be solved within minutes. The main contributions of this paper are summarized as follows.

1) A systematic framework that properly combines the advantages of the analytical function and the NSTD is proposed. In particular, the NSTD is used to adjust the trial UC solution given by the analytical function in such a way that the adjusted UC solution is likely to satisfy all the constraints.

2) A corrective method for making the adjusted UC solution completely feasible is proposed based on a constructed feasibility-testing problem. Specifically, the feasibility-testing problem is constructed in such a way that it yields information about the gradient of the Lagrangian dual function, and this information is essential for constructing a feasible update of the trial UC solution. In addition, the updated feasible solution is of high quality in terms of the generation cost.

The rest of this paper is organized as follows. Section II describes the basic UC model and the fast trial UC solution in [15]. Section III firstly shows the main idea of the method and then elaborates on it. This method consists of fast settings of the trial UC solution, an adjustment method for the trial UC solution, and a method for constructing a feasible sub-optimal UC solution. Section IV presents thorough numerical results. Section V concludes this paper.

## II. BASIC UC PROBLEM AND A FAST TRIAL UC SOLUTION

In this section, the basic formulation of the UC problem and the fast trial UC solution based on LR are introduced for later analysis.

The security-constrained UC model is adopted in this paper and is presented as follows.

$$\min_{z,p} \sum_i S_i(z_i) + \sum_{i,t} C_i(z_{i,t}, p_{i,t}) \quad (1)$$

$$\text{s.t.} \quad \sum_i p_{i,t} - \sum_m d_{m,t} = 0, \forall t; \quad (2)$$

$$\sum_i \Gamma^U_{l,i} p_{i,t} - \sum_m \Gamma^D_{l,m} d_{m,t} \leq F_l, \forall l,t; \quad (3)$$

$$-\sum_i \Gamma^U_{l,i} p_{i,t} + \sum_m \Gamma^D_{l,m} d_{m,t} \leq F_l, \forall l,t; \quad (4)$$

$$z_{i,t} \underline{P}_i \leq p_{i,t} \leq z_{i,t} \overline{P}_i, \forall i,t; \quad (5)$$

$$\Delta^-_i(z_{i,t}, z_{i,t-1}) \leq p_{i,t} - p_{i,t-1} \leq \Delta^+_i(z_{i,t}, z_{i,t-1}), \forall i,t; \quad (6)$$

$$Uz \leq R, z \in \{0,1\}^{I \times T} \quad (7)$$

The objective (1) is to minimize the total generation cost, including the start-up costs and the fuel costs. Equation (2) is the power balance constraint. Inequations (3)-(4) are the transmission capacity/security constraints based on the direct current (DC) power flow model. Constraints (5) and (6) represent the generation capacities and ramp up/down limits. Equation (7) is the minimum up/down time constraints. According to [29], spinning reserve constraints can also be incorporated in (7). The system-wide constraints (2)-(4) and the unit-wise constraints (5)-(7) incorporate representative coupling relationships of the variables. Detailed information about this UC model could be referred to in [3].

For (1)-(7), the LR methods usually relax the system-wide constraints, and the corresponding relaxed problem is as (8), where the multipliers $\lambda_{0,t}$, $\lambda^+_{l,t}$, and $\lambda^-_{l,t}$ are associated with the constraints (2), (3), and (4), respectively. Then, the LR methods generate a series of trial solutions that converge to the optimal solution by alternately solving (8) and maximizing $q(\lambda)$ over $\lambda$.

However, the LR methods yield many computationally demanding unit sub-problems for large-scale UC problems [16].

$$\begin{cases} q(\lambda) = \min_{z,p} \sum_i S_i(z_i) + \sum_{i,t} C_i(z_{i,t}, p_{i,t}) + \\ \quad \sum_t \lambda_{0,t}(\sum_i p_{i,t} - \sum_m d_{m,t}) + \\ \quad \sum_{l,t} \lambda^+_{l,t}(\sum_i \Gamma^U_{l,i} p_{i,t} - \sum_m \Gamma^D_{l,m} d_{m,t} - F_l) + \\ \quad \sum_{l,t} \lambda^-_{l,t}(-\sum_i \Gamma^U_{l,i} p_{i,t} + \sum_m \Gamma^D_{l,m} d_{m,t} - F_l) \\ \text{s.t. } \Delta^-_i(z_{i,t}, z_{i,t-1}) \leq p_{i,t} - p_{i,t-1} \leq \Delta^+_i(z_{i,t}, z_{i,t-1}), \forall i,t; \\ \quad z_{i,t} \underline{P}_i \leq p_{i,t} \leq z_{i,t} \overline{P}_i, \forall i,t; \\ \quad Uz \leq R, z \in \{0,1\}^{I \times T} \end{cases} \quad (8)$$

This computational difficulty can be well addressed by further relaxing the unit-wise constraints (6)-(7) [15]. In particular, after the system-wide constraints are relaxed and the unit-wise constraints are ignored, the Lagrange dual function for (1)-(7) is as follows.

$$\hat{q}(\lambda) = \min_{z,p} \sum_i S_i(z_i) + \sum_{i,t} C_i(z_{i,t}, p_{i,t}) + \\ \sum_t \lambda_{0,t}(\sum_i p_{i,t} - \sum_m d_{m,t}) + \\ \sum_{l,t} \lambda^+_{l,t}(\sum_i \Gamma^U_{l,i} p_{i,t} - \sum_m \Gamma^D_{l,m} d_{m,t} - F_l) + \\ \sum_{l,t} \lambda^-_{l,t}(-\sum_i \Gamma^U_{l,i} p_{i,t} + \sum_m \Gamma^D_{l,m} d_{m,t} - F_l) \quad (9)$$

$$\text{s.t.} \quad z_{i,t} \underline{P}_i \leq p_{i,t} \leq z_{i,t} \overline{P}_i, z_{i,t} \in \{0,1\}, \forall i,t \quad (10)$$

where the hat sign "^" is for denoting $\hat{q}(\lambda)$ as a relaxation of $q(\lambda)$. Then, given the optimal solution to (9)-(10), the multiplier vector $\lambda$ (i.e., the dual variables) is to be optimized with the following dual problem.

$$\max_\lambda \hat{q}(\lambda) \quad (11)$$

$$\text{s.t.} \quad \lambda_0 \text{ free}, \lambda_l \geq 0, \forall l \quad (12)$$

At this time, the optimal solution pair $(\hat{z}^*, \hat{p}^*, \hat{\lambda}^*)$ to (9)-(12) can be easily obtained. Specifically, the problem (9)-(10) is separable both in periods and in units, and this problem can be separated into $I \times T$ unit-period sub-problems. The unit-period sub-problem of unit $i$ at period $t$ is as follows.

$$\hat{q}_{i,t}(\lambda_t) = \min_{z_{i,t}, p_{i,t}} C_i(z_{i,t}, p_{i,t}) + \beta_{i,t}(\lambda_t) p_{i,t} \quad (13)$$

$$\text{s.t.} \quad z_{i,t} \underline{P}_i \leq p_{i,t} \leq z_{i,t} \overline{P}_i, z_{i,t} \in \{0,1\} \quad (14)$$

where $\beta_{i,t}(\lambda_t) = -\lambda_{0,t} + \sum_l (\lambda^+_{l,t} - \lambda^-_{l,t}) \Gamma^U_{l,i}$ (of which the constant part is omitted). With any given multiplier vector $\lambda_t$, the associated optimal binary solution $\hat{z}^*_{i,t}(\lambda_t)$ to (13)-(14) can be expressed as the following analytical function, and the dispatch solution $\hat{p}^*_{i,t}(\lambda_t)$ of this problem can also be easily obtained [15].

$$\hat{z}^*_{i,t}(\lambda_t) = \begin{cases} 1, & \text{if } \beta_{i,t}(\lambda_t) < \beta^0_i \\ 0 & \text{if } \beta_{i,t}(\lambda_t) \geq \beta^0_i \end{cases} \quad (15)$$

where the threshold $\beta^0_i$ for determining the binary solution $\hat{z}^*_{i,t}(\lambda_t)$ is (16). It needs to be clarified that, (13) does not include the commitment cost since it becomes constant in the deduction of (15) [15].

$$\beta^0_i = \max\{-a_i - b_i / \overline{P}_i, -a_i - b_i / \underline{P}_i\} \quad (16)$$

With the solution pair $(\hat{z}^*, \hat{p}^*)$, the dual problem (11)-(12) becomes the following single-level optimization problem.

$$\{max_\lambda \hat{q}(\lambda) : \lambda_{0,t} \text{ free}, \lambda^+_{l,t}, \lambda^-_{l,t} \geq 0, \forall l,t\} \quad (17)$$

which is separable in periods and whose optimal solution $\hat{\boldsymbol{\lambda}}^*$ can be easily obtained with $T$ linear programs [15]. Then, with this optimal multiplier vector $\hat{\boldsymbol{\lambda}}^*$, the analytical function (15) directly gives the fast trial UC solution $\hat{z}^*(\hat{\boldsymbol{\lambda}}^*)$ to (9)-(10). This solution fully considers the system-wide constraints; however, due to the relaxation, it is generally infeasible and needs adjustments.

## III. THE PROPOSED FAST SOLUTION METHOD FOR THE LARGE-SCALE UC PROBLEM

### A. Main Idea of the Method

The proposed method aims at providing a fast sub-optimal UC solution to the basic UC problem based on efficiently and properly adjusting its trial UC solutions and the associated multipliers, iteratively (see Fig. 1).

The method is explained in the following steps.

**Step 0. Initialize.** Set $k=0$, set the initial LR multiplier vector $\hat{\boldsymbol{\lambda}}_0^*$ as $\hat{\boldsymbol{\lambda}}^*$; go to step 1.

**Step 1. Calculate trial UC.** With $\hat{\boldsymbol{\lambda}}_k^*$, calculate the trial UC solution (denoted as $\hat{z}(\hat{\boldsymbol{\lambda}}_k^*)$) with (15). According to (9)-(10), this trial UC fully considers the system-wide constraints, but it may be infeasible to the basic UC problem and will be adjusted. Go to step 2.

**Step 2. Adjust the trial UC** $\hat{z}(\hat{\boldsymbol{\lambda}}_k^*)$. With NSTD, the trial UC $\hat{z}(\hat{\boldsymbol{\lambda}}_k^*)$ is adjusted to a new UC solution (denoted as $\bar{z}(\hat{\boldsymbol{\lambda}}_k^*)$) in such a way that, this new UC $\bar{z}(\hat{\boldsymbol{\lambda}}_k^*)$ is feasible to the unit-wise constraints and has the minimum difference with the trial UC. At this time, this new UC $\bar{z}(\hat{\boldsymbol{\lambda}}_k^*)$ is likely to be feasible to all the constraints. Go to step 3.

**Step 3. Test and improve the feasibility of the new UC.** Firstly, set the binary variables in the basic UC problem (1)-(7) as the new UC $\bar{z}(\hat{\boldsymbol{\lambda}}_k^*)$. This yields a linear program that may not be feasible. Secondly, add non-negative slack variables to this linear program to replace the possible violations of its system-wide constraints, and set its minimizing objective as the sum of the slack variables. To this end, this linear program must be feasible and serves as a feasibility-testing problem for the new UC. Thirdly, if the slack variables (violations) are all minimized to zeros or their sum is less than a convergence tolerance, the new UC is regarded as feasible and the iteration stops; otherwise, these minimized violations turn out to constitute an approximate sub-derivative $\boldsymbol{g}_{\hat{\lambda}_k^*}$ of the Lagrangian dual function $q(\boldsymbol{\lambda})$ at the multiplier vector $\hat{\boldsymbol{\lambda}}_k^*$. Go to step 4.

**Step 4. Update the multipliers.** Update $\hat{\boldsymbol{\lambda}}_k^*$ to $\hat{\boldsymbol{\lambda}}_{k+1}^*$ with the sub-derivative $\boldsymbol{g}_{\hat{\lambda}_k^*}$ and with a properly optimized step size; set $k=k+1$, and go to step 1.

### B. Adjustment Method for the Fast Trial UC Solution

In Fig. 1, the trial UC is already provided by (15); this subsection explains how to adjust the trial UC to the new UC. Namely, the new UC is obtained with the NSTD in such a way that it strictly satisfies the unit-wise constraints and has the minimum difference from the trial UC.

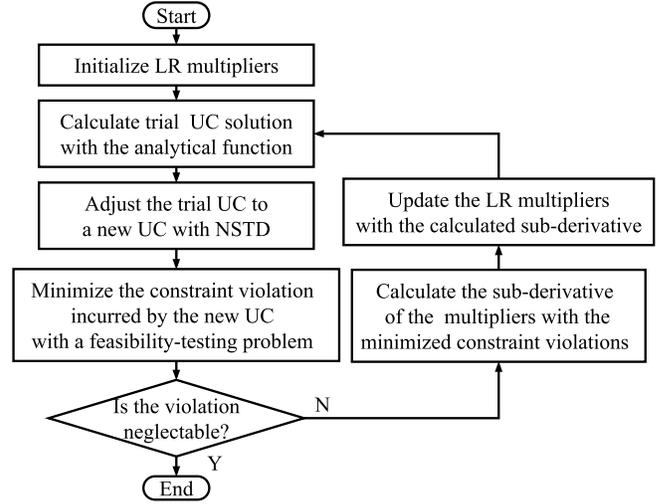

Fig. 1. The basic idea of the proposed method.

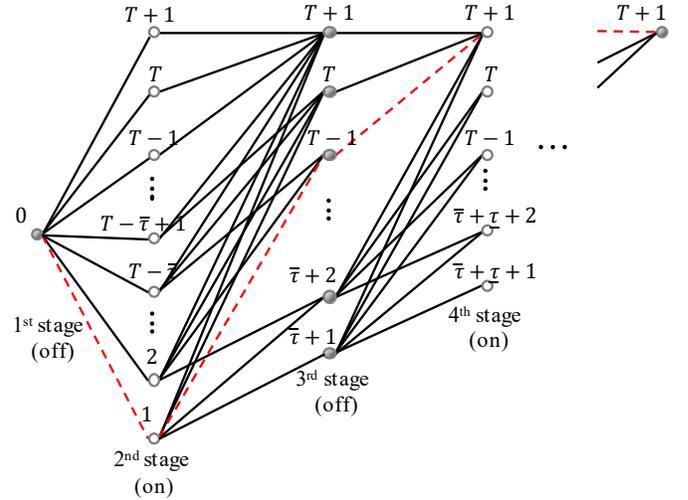

Fig. 2. The new state transition diagram of the unit (the unit has been off for the minimum down time before the scheduling horizon).

The NSTD serves as the basis for the adjustment. As seen in Fig. 2, this state transition diagram for each unit consists of hundreds of transition states scattered in multiple stages. Each of the states represents a scheduling period and is labeled with the corresponding serial number of the period; each stage corresponds to a decision of turning the unit on (or off) for a time span (which is longer than the minimum up/down time), and each of the paths connecting the two states labeled with 0 and T+1 yields a UC decision chain for the unit during the horizon. For example, the path of the red dotted line segments in Fig. 2 corresponds to a UC decision chain that starts from turning the unit on at period 1; then, this "on" decision lasts for T-1 hours ($\geq \bar{\tau}_i = 4$) until the period T-1 at the $3^{rd}$ stage; then, the next "off" decision lasts until the period T+1 at the $4^{th}$ stage (which is in the next scheduling horizon). It is pointed out that the on/off status of the unit determined by each of the paths satisfies the unit-wise constraints [17]. And since the trial UC may violate the unit-wise constraints, there may be a difference between the on/off status of each path and the trial UC. These differences are essential in adjusting the trial UC to the new UC.

The trial UC is to be adjusted to the on/off status determined by a particular one of the paths; also, this on/off status has the minimum difference from the trial UC, and is actually the new UC. This is realized by applying dynamic programming on NSTD with the following procedures.

Firstly, to apply dynamic programming on NSTD, the transition cost associated with each of the edges is defined. In particular, the on/off status determined by any one of the edges may also be different from the trial UC of the unit at the corresponding periods. This difference is defined as the transition cost of this edge, as follows.

$$\Omega_{i,e}(\hat{z}^*(\lambda), z^e) = \sum_{t_1}^{t_2-1} \left| \hat{z}^*_{i,t}(\lambda) - z^e_{i,t} \right| \tag{18}$$

where $t_1, t_2$ are endpoint periods associated with the edge $e$, and $\{z^e_{i,t} \mid t_1 \leq t \leq t_2 - 1, t \text{ integer}\}$ is the on/off status determined by this edge.

Secondly, the transition cost associated with any one of the paths $c$ in $G_i$ is defined as (19), which is the sum of the transition costs (18) of the edges therein.

$$\Omega_{i,c}(\hat{z}^*(\lambda), z^c) = \sum_e \Omega_{i,e}(\hat{z}^*(\lambda), z^e), e = 1, 2, \ldots \tag{19}$$

where $e$ is the index of the edges in the path $c$. Then, the optimal path $c^*$ with the minimum total cost $\Omega_{i,c^*}$ is obtained by solving the following problem with dynamic programming.

$$c^* = \arg\min_{c \in G_i} \Omega_{i,c}(\hat{z}^*(\lambda), z^e) \tag{20}$$

It is seen that, the on/off status determined by the optimal path $c^*$ does have the minimum difference with the trial UC $\hat{z}^*(\hat{\lambda}^*)$, and it is defined as the aforementioned new UC $\bar{z}^*(\hat{\lambda}^*)$ of the unit. In this way, this new UC fully considers the system-wide constraints besides strictly satisfying the unit-wise constraints, and it is of high quality in terms of the fuel cost. In addition, the commitment cost omitted in (13)-(14) can be easily calculated for this new UC at this time.

Although the new UC is likely to satisfy all the constraints, it does not guarantee to satisfy the system-wide constraints. Therefore, further adjustments to the new UC are needed.

*C. Construction Method for Feasible UC Solution*

If the new UC is infeasible to the basic UC problem, it can only have violations of the system-wide constraints. In this case, the new UC is to be adjusted in this subsection to accommodate these violations. This is realized by adjusting the multipliers that determine the new UC $\bar{z}^*(\hat{\lambda}^*)$.

The feasibility of the new UC can be efficiently examined with the feasibility-testing problem. This problem is a linear program derived by fixing the binary variables in the basic UC problem as the new UC. If the new UC is feasible, in the feasible region of this linear program, there exists a dispatch solution point that has zero violations of all the constraints, including the system-wide ones; otherwise, there is no such point, and violations of the system-wide constraints are inevitable. In the latter case, the minimum total violation of the system-wide constraints is to be obtained with the feasibility-testing problem, as follows.

$$\tilde{q}(\lambda) = \min_{p,v} \sum_{l,t}(v^+_{l,t} + v^+_{\bar{l},t}) + \sum_t(v^+_{0,t} + v^+_{\bar{0},t}) \tag{21}$$

$$\text{s.t. } \sum_i \Gamma^U_{l,i} p_{i,t} - \sum_m \Gamma^D_{l,m} d_{m,t} - v^+_{l,t} \leq F_l, \forall l,t; \tag{22}$$

$$-\sum_i \Gamma^U_{l,i} p_{i,t} + \sum_m \Gamma^D_{l,m} d_{m,t} - v^+_{\bar{l},t} \leq F_l, \forall l,t; \tag{23}$$

$$\sum_i p_{i,t} - \sum_m d_{m,t} + v^+_{0,t} - v^+_{\bar{0},t} = 0, \forall t; \tag{24}$$

$$\bar{z}^*_{i,t}(\hat{\lambda}^*_k) \underline{P}_i \leq p_{i,t} \leq \bar{z}^*_{i,t}(\hat{\lambda}^*_k) \bar{P}_i, \forall i,t; \tag{25}$$

$$\Delta^-_i(\bar{z}^*_{i,t}, \bar{z}^*_{i,t-1}) \leq p_{i,t} - p_{i,t-1} \leq \Delta^+_i(\bar{z}^*_{i,t}, \bar{z}^*_{i,t-1}), \forall i,t; \tag{26}$$

$$v^+_{0,t}, v^+_{\bar{0},t}, v^+_{l,t}, v^+_{\bar{l},t} \geq 0, \forall l,t; \tag{27}$$

where the nonnegative auxiliary vector variable $v$ is used to represent the non-negative violations (residuals) of the system-wide constraints. The model (21)-(27) must be feasible if the basic UC problem is feasible. Then, as in the objective (21), all the auxiliary variables are minimized to see if there exists a dispatch solution $\bar{p}^*(\hat{\lambda}^*_k)$ under $\bar{z}^*(\hat{\lambda}^*_k)$ that has zero total violation of the system-wide constraints. Once the optimal objective is zero, the optimal solution pair $(\bar{z}^*(\hat{\lambda}^*_k), \bar{p}^*(\hat{\lambda}^*_k))$ of the above feasibility-testing problem must be a feasible sub-optimal solution to the basic UC problem.

If the optimal objective is greater than zero, the optimized violations (i.e., the optimal values of the slack variables) are to be used to adjust the multipliers and thus the trial UC solutions, iteratively. According to LR multiplier theory, the non-zero optimal values of the auxiliary variables $v^*$ turn out to constitute an approximate sub-derivative of $q(\lambda)$ at the multiplier vector $\hat{\lambda}^*_k$, as follows.

$$g_{\hat{\lambda}^*_k} = (v^{+*}_{0,t} - v^{+*}_{\bar{0},t}, v^{+*}_{l,t}, v^{+*}_{\bar{l},t}) \tag{28}$$

where $k$ is the index of the iterations for the adjustments. This approximate sub-derivative is essential in updating the multiplier vector $\hat{\lambda}^*_k$ to a new multiplier vector $\hat{\lambda}^*_{k+1}$, with the sub-gradient method as follows.

$$\hat{\lambda}^*_{k+1} = \hat{\lambda}^*_k + c_k g_{\hat{\lambda}^*_k} \tag{29}$$

where $c_k$ is a step size specified as $c_0 / k > 0$ ($c_0$ is configured with the method in the appendix). Compared with $\hat{\lambda}^*_k$, $\hat{\lambda}^*_{k+1}$ derives a UC solution achieving a better objective of the dual problem (8). Namely, this UC solution has a decreased violation of the system-wide constraints, with respect to the generation cost. This also holds for the trial UC $\hat{z}(\hat{\lambda}^*_{k+1})$ and the new UC $\bar{z}^*(\hat{\lambda}^*_{k+1})$ calculated with the multipliers $\hat{\lambda}^*_{k+1}$, (15), and (20). In other words, besides attaining high quality in terms of the generation cost, the new UC $\bar{z}(\hat{\lambda}^*_{k+1})$ has improved feasibility compared with that of the previous new UC $\bar{z}(\hat{\lambda}^*_k)$.

A feasible new UC is obtained with iterations of the above adjustments. Under the framework of sub-gradient methods, multiple iterations of the adjustments may be required to derive a feasible new UC. Specifically, in the possible case that the new UC $\bar{z}(\hat{\lambda}^*_{k+1})$ is infeasible, the feasibility-testing problem will again derive an approximate sub-derivative $g_{\hat{\lambda}^*_{k+1}}$. This approximate sub-derivative $g_{\hat{\lambda}^*_{k+1}}$ is then used to update $\hat{\lambda}^*_{k+1}$ to $\hat{\lambda}_{k+2}^*$ as the way of updating $\lambda_k^*$ to $\lambda_{k+1}^*$. With (15), (18)-(20), this newest multiplier vector $\hat{\lambda}_{k+2}^*$ will again yield a new UC $\bar{z}(\hat{\lambda}_{k+2}^*)$ with further improved feasibility. The above process illustrated in Fig. 1 repeats until a feasible new UC appears, or the following total violation of the system-wide constraints is less than a specified convergence tolerance $\varepsilon$.

$$V = \sum_{l,t}(v_{l,t}^{+*}+v_{l,t}^{-*}) + \sum_{t}(v_{0,t}^{+*}+v_{0,t}^{-*}) \leq \varepsilon \qquad (30)$$

It is verified in numerical tests that only after a few iterations, the solution $(\bar{z}^*, \bar{p}^*)$ from (21)-(27) satisfies all the constraints in the basic UC problem, or (30) is strictly satisfied. From now on, the proposed fast solution method based on DP and LR is denoted as DPLR.

## IV. NUMERICAL TESTS

In this section, an overview of the performance of the method is described first. Then, the solution time and the solution quality are tested under normal constraints. Finally, the method is tested under strengthened constraints. All these numerical tests are implemented with Matlab 2018b and the solver Gurobi 9.5.0 running on a desktop with an Intel 3.3GHz CPU and 32 GB RAM.

### A. Basic Information and Performance Summary

The method is tested on cases of the IEEE 24-bus system, 118-bus system, a modified Polish 2383-bus system, and a modified French 6468-bus system, with varied tightness of the coupling constraints (i.e., different levels of loads, ramp limits, minimum up/down time, and transmission limits). The performance is also benchmarked against that of the Gurobi solver. Details of the tested cases are available in [30]-[33]. Basic information about these cases and the parameters of DPLR are listed in Table I and Table II, respectively. In Table II, the column "$\Delta^{\pm}/\bar{P}_i$" is the ratio of the ramp limit to the corresponding maximum generation capacity of the unit; the column "$\bar{\tau}_i$" shows the maximum of the minimum up/down times of the unit; $\varepsilon$ is the specified convergence tolerance of DPLR; $\delta, u$ are given quantities for optimizing the initial step size in the appendix.

With the convergence tolerance $\varepsilon$ in (30) set as $10^{-4}$, the overall solution time of DPLR is on average less than one-tenth of the time consumed by the solver. In particular, 1) it is found that the computation time of DPLR for solving problems of the Polish 2383-bus case is within several minutes; 2) DPLR only requires a few iterations to provide a feasible sub-optimal solution even if the test cases are with very tight constraints; 3) the solution of DPLR achieves competitive generation cost compared with that of the solver.

### B. Analysis of the Computational Efficiency

In this subsection, the computational efficiency is tested by examining the number of iterations and the solution time of DPLR and by comparing these results with that of the solver.

Each of the iterations in DPLR is fast due to its computationally simple operations. The feasibility-testing problem consumes the major computation time of DPLR; however, it is very easy to solve since it only has much fewer continuous variables compared with the mixed variables in the basic UC problem. In addition, redundant constraints of the feasibility-testing problem are further eliminated to reduce the computational complexity, without changing its optimal solution [34]. This benefits the efficiency of DPLR. As in Table III, the time consumed by one of the iterations is in the column "Time", of which the sub-columns "ERC" and "No ERC" show the solution times of DPLR with and without eliminating the redundant constraints, respectively; and the time consumed by the Gurobi solver is in the sub-column "Gurobi", in which the "--" sign means that the solver runs out the RAM. The sub-columns "Var" and "Constr" in column "FT" list the numbers of variables and the numbers of constraints in the feasibility-testing problem, and the value before a slash sign in the column "Constr" is the number of the constraints after the elimination, and value after a slash sign is the number of the full constraints. The sub-columns "Var" and "Constr" in the column "BUC" list the number of the mixed variables and the constraints in the basic UC problem.

TABLE I
BASIC INFORMATION

| Sys. | $N$ | $L$ | $T$ | $\sum_i \bar{P}_i$(MW) | $\sum_{k,t} d_{k,t}$(MW) |
|---|---|---|---|---|---|
| 24-bus | 32 | 38 | 24 | 3405 | 2850 |
| 118-bus | 54 | 186 | 24 | 13373 | 6000 |
| 2383-bus | 327 | 2896 | 24 | 29934 | 24558 |
| 6468-bus | 400 | 9000 | 24 | 115596 | 52535 |

TABLE II
PARAMETERS OF THE CASES

| $d_{k,t}$ (MW) | $\Delta^{\pm}/\bar{P}_i$ | $F_l$ (MW) | $\bar{\tau}_i$ (h) | $\varepsilon$ | $\delta$ | $u$ |
|---|---|---|---|---|---|---|
| 100%Peak Load | $\leq 0.6$ | 100% | $\leq 6$ | $10^{-4}$ | 0.5 | 1 |

TABLE III
COMPARISON BETWEEN THE SOLVER AND DPLR

| Sys. | BUC | | FT | | Time(sec.) | | |
|---|---|---|---|---|---|---|---|
| | Var | Constr | Var | Constr | ERC | No ERC | Gurobi |
| 24 | 1536 | 6456 | 577 | 1660/3002 | 0.45 | 0.61 | 0.8 |
| 118 | 2592 | 16728 | 910 | 1892/10772 | 0.94 | 1.24 | 1.0 |
| 2383 | 15696 | 186120 | 2062 | 14619/143190 | 7.15 | 18.09 | 2463.5 |
| 6468 | 19200 | 486120 | 2552 | 21082/442776 | 15.50 | 57.52 | -- |

It is seen from Table III that: 1) the feasibility-testing problem has much fewer variables and constraints compared to the basic UC problem; 2) in addition, elimination of the redundant constraints distinctively improves the efficiency; 3) thus, it is seen that the computation time of one of the iterations is very small and grows almost linearly with the system size. This indicates that the DPLR can be used to solve the problem of even larger systems. These advantages of DPLR are attributed to the fact that all its computations are simple operations with low computational complexity.

Then, the number of iterations that DPLR requires to solve the problem is obtained under different settings of the initial step size. This is because the initial step size of sub-gradient methods has a significant influence on the number of iterations, and it usually needs to be tuned. Nevertheless, with the method in the appendix, the initial step size $c_0$ can be properly optimized as $c_0^*$ without experimentally tuning it. To test both the number of the required iterations and the effectiveness of the optimized initial step size $c_0^*$, 30 exponentially increasing initial step sizes are generated within the range of $[0,c_0^*]\bigcup(c_0^*,1]$, and the number of the total iterations of DPLR at each generated step size is obtained. As in Fig. 3, the horizontal axis labels the logarithm of the step size, and the vertical axis shows the number of iterations; arrows are added to indicate the optimized step size $c_0^*$ of each system. The DPLR is terminated if the number of the iterations exceeds 20; and the markers in the legend of Fig. 3 label the step sizes at which a feasible DPLR solution is obtained within 20 iterations.

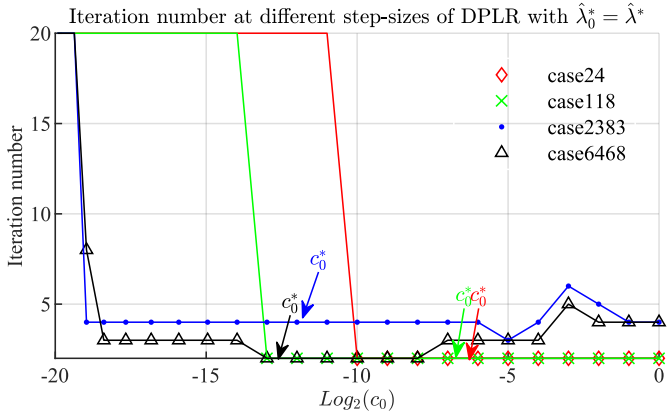

Fig. 3. Convergence speed.

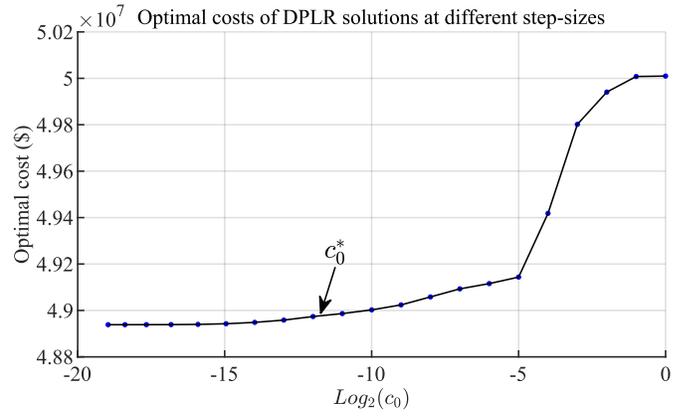

Fig. 5. Case2383: optimal costs of DPLR solutions at different step sizes.

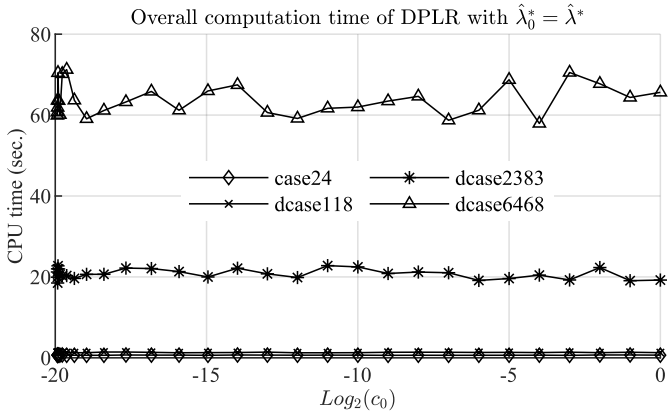

Fig. 4. Overall solution time of DPLR.

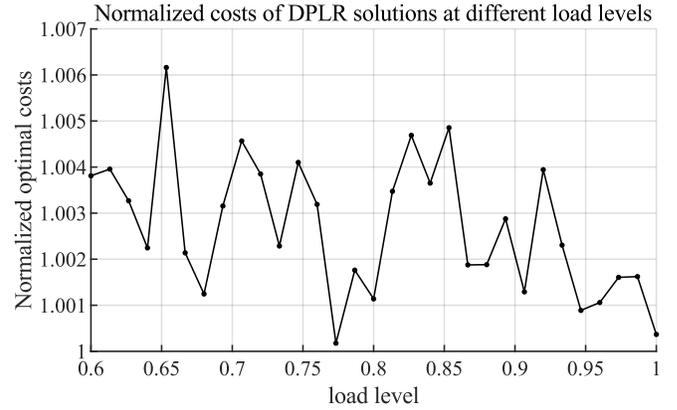

Fig. 6. Case2383: normalized costs DPLR solutions at different load levels.

It is seen from Fig. 3 that: 1) the typical number of the total iterations is 2 to 4. This combined with TABLE III suggests that DPLR is fast in providing a feasible solution; 2) for cases of smaller systems, the initial step size inducing a feasible DPLR solution is typically bigger. This is mainly because, for these cases the non-negative violations (sub-gradient) are small, so, progressive updates of the multipliers can only be obtained with a relatively bigger step size; 3) at the initial step size around the optimized step size $c_0^*$, the iteration number is all smaller than 4; this verifies the effectiveness of $c_0^*$ that $c_0^*$ can serve as a proper initial step size. 4) at the initial step size smaller than the optimized initial step size $c_0^*$, DPLR usually requires more iterations. This is mainly because this step size is too small to make DPLR has progressive updates of the trial solutions. 5) at the initial step size that is much bigger than the optimized step size $c_0^*$, DPLR still requires only a few iterations. However, it is found that DPLR with a much bigger step size may yield a solution with more generation cost, which will be discussed in the next subsection.

In Fig. 4, the overall solution time of DPLR at each previously generated step size is illustrated. It is seen that:1) for the small cases, the computation time of DPLR is comparable with that of the solver. 2) for large-scale cases, the computation time of DPLR is on average one-tenth of the time consumed by the solver. This is expectable since only large-scale cases can serve to distinctly examine the advantages of DPLR.

### C. Analysis of the Solution Quality

In practice, the generation cost of a scheduling plan must be considered by the system operators. This aspect of DPLR is discussed in this subsection.

The generation cost of a UC solution is accordingly defined as follows.

$$\mathbf{\Lambda}_u(z, p) = \sum_i S_i(z_i) + \sum_{i,t} C_i(z_{i,t}, p_{i,t}) \qquad (31)$$

And for convenience, we also define the following normalized cost of the DPLR solution.

$$\mathbf{N}_c = \mathbf{\Lambda}_c(\bar{z}^*, \bar{p}^*) / \mathbf{\Lambda}_c(z^*, p^*) \qquad (32)$$

where $(\bar{z}^*, \bar{p}^*)$ is a DPLR solution and $(z^*, p^*)$ is a solver solution.

At the initial step size bigger than the optimized step size $c_0^*$, DPLR solution incurs more generation cost than the cost of the DPLR solution associated with $c_0^*$. Specifically, at each previously generated initial step size that derives feasible solutions for the 2383-bus system, we solve the UC problem with DPLR and obtain the corresponding optimal cost of (31). Then, the optimal cost is illustrated in Fig. 5.

Fig. 5 shows that: 1) a bigger initial step size incurs a heavier generation cost of DPLR solution. This is attributed to the fact that the bigger the step size is, the sharper the feasibility of the trial solution is improved, and the more the generating resource is committed. 2) on the other side, the costs associated with the step sizes that are smaller than $c_0^*$ are almost the same as the cost associated with the optimized step size $c_0^*$. Tests of other systems show similar results as above.

The cost of DPLR solution is then compared with that of the solver at different load levels. For the 2383-bus cases with the parameters in Table II and with the logarithmic initial step size being $-13.3$, the curve in Fig.6 shows the normalized optimal cost $\mathbf{N}_c$ of DPLR solution at a load level ranging from 60% to 100% of the peak load.

It is seen from Fig. 6 that, the optimal DPLR solution has a slightly increased generation cost (less than 1%) compared to that of a solver solution; this indicates that the DPLR solution is of high quality in terms of generation costs. Tests of other systems show similar results.

### D. Performance under Tighter Constraints

This subsection shows that DPLR still has almost the same performance as above even with much tighter constraints, such as heavier loads, lower transmission limits, smaller ramp rates, and longer minimum up/down time.

The tightness of the load, transmission limits, ramp limits, and minimum up/down time are measured by the following four scaling factors.

$$s_d = d^\uparrow / d, s_M = \overline{\tau}^\uparrow / \overline{\tau}, s_F = F^\downarrow / F, s_R = \Delta^{\pm\downarrow} / \Delta^\pm \quad (33)$$

where $d^\uparrow$ and $\overline{\tau}^\uparrow$ represent the uniformly increased loads and minimum up/down time, and $F^\downarrow$, $\Delta^{\pm\downarrow}$ represent the uniformly decreased transmission limits and ramp up/down limits, respectively. With each kind of constraint becoming tighter, there will be an extreme value of the corresponding scaling factor at which the DPLR (or the solver) will fail to generate a feasible solution; such extreme values of the scaling factors are denoted as $\overline{s}_d, \overline{s}_M, \underline{s}_F, \underline{s}_R$, respectively. As in Table IV, the values before the slash sign are the extreme scaling factors of the solver, the remaining values are the extreme scaling factors of DPLR. Notice that there is no value for the solver with the 6468-bus cases since under which the RAM runs out, and the value in the column $\overline{s}_M$ is truncated if the corresponding scaled minimum up/down time exceeds the length of the scheduling horizon.

It is seen from Table IV that, 1) the extreme values of the solver are slightly "tighter" than that of DPLR. This means DPLR cannot handle cases where the feasible region is extremely small; 2) however, tests show that such extreme cases also make the solver consume excessive computing time and resources, and these cases are rare in practice.

Then, after scaling the transmission limits (other constraints remain normal) of the systems with the corresponding extreme values $\underline{s}_F$ in Table IV, there are many versions of the problem that have varied tightness of the constraints. These problems are solved by DPLR, and the associated total number of iterations and the normalized optimal costs $\mathbf{N}_c$ of DPLR are listed in the columns "$k_a$" and "cost" in Table V. Table VI shows corresponding results of DPLR after scaling the ramp limits with $\underline{s}_R$, and Table VII shows the corresponding results of DPLR after scaling the minimum up/down time with $\overline{s}_M$. In Table V, VI, and VII, the second column lists the corresponding extreme scaling factors $\overline{s}_M, \underline{s}_F$, and $\underline{s}_R$; the "Y" in the "solved?" columns means that feasible sub-optimal solutions are obtained within 20 iterations, and "N" means otherwise. The solution time of DPLR is in the last column.

TABLE IV
LIMITS OF THE SCALING FACTORS

|  | $\overline{s}_d$ | $\overline{s}_M$ | $\underline{s}_F$ | $\underline{s}_R$ |
| --- | --- | --- | --- | --- |
| 24-bus | 1.2/1.1 | 1.7/1.3 | 0.4/0.4 | 0.3/0.4 |
| 118-bus | 2.2/2.0 | 4.6/4.3 | 0.4/0.8 | 0.6/0.6 |
| 2383-bus | 1.1/1.0 | 5.0/5.0 | 0.8/0.8 | 0.5/0.6 |
| 6468-bus | --/1.1 | --/3.0 | --/0.9 | --/0.6 |

TABLE V
PERFORMANCE OF DPLR WITH TIGHTER TRANSMISSION LIMITS

| System | $\underline{s}_F$ | Solved? | $k_a$ | error | CPU time (sec) |
| --- | --- | --- | --- | --- | --- |
| 24bus | 0.4 | Y | 2 | 1.01 | 2.0 |
| 118bus | 0.8 | Y | 3 | 1.007 | 5.9 |
| 2383bus | 0.8 | Y | 5 | 1.004 | 59.3 |
| 6468bus | 0.9 | Y | 2 | 1.006 | 109.6 |

TABLE VI
PERFORMANCE OF DPLR WITH TIGHTER RAMP LIMITS

| System | $\underline{s}_R$ | Solved? | $k_a$ | error | CPU time (sec) |
| --- | --- | --- | --- | --- | --- |
| 24bus | 0.4 | Y | 2 | 1.003 | 2.0 |
| 118bus | 0.6 | Y | 2 | 1.001 | 4.5 |
| 2383bus | 0.6 | Y | 6 | 1.006 | 67.5 |
| 6468bus | 0.6 | Y | 2 | 1.009 | 107.7 |

TABLE VII
PERFORMANCE OF DPLR WITH TIGHTER MINIMUM UP/DOWN TIME

| System | $\overline{s}_M$ | Solved? | $k_a$ | error | CPU time (sec) |
| --- | --- | --- | --- | --- | --- |
| 24bus | 1.3 | N | 20 | 1.01 | -- |
| 118bus | 4.3 | Y | 3 | 1.02 | 6.9 |
| 2383bus | 5.0 | Y | 8 | 1.008 | 66.1 |
| 6468bus | 3.0 | Y | 2 | 1.005 | 105.6 |

### V. CONCLUSION

Unit commitment problems for large-scale power systems are difficult to solve mainly due to their non-convex and non-differentiable nature. This difficulty in LR-based methods lies in the feasibility issues and the computational problems of their trial solutions. The computational difficulty can be greatly alleviated by the analytical function that gives a fast trial UC solution under given multipliers. With this analytical function, a fast solution method for large-scale UC problems is established in this paper based on Lagrange relaxation and dynamic programming. The method searches for a sub-optimal UC solution with the NSTD rather than the tree structure used in MILP methods; and by constructing proper directions for updating the trial solutions, the method usually requires only a few iterations to provide the solution.

The numerical results show that DPLR is on average ten times faster than the commercial solver in solving large-scale UC problems, and DPLR solution is of high quality in terms of the generation cost. The generalization of DPLR for solving UC problems with uncertainty is in progress. DPLR is inherently generating a fast solution to large-scale UC problems under a single scenario of loads. In the future, the feasibility of DPLR solution under many of the scenarios or even all the scenarios of the uncertain loads will be assured.

## APPENDIX

For sub-gradient methods, the initial step size usually needs to be tuned for a tested system, and this tuned step size may not be applicable to other systems. This results in re-tuning of the initial step size once the system is different. To find a proper initial step size without experimentally tuning it, an ideal approach is to start DPLR from a good-enough trial solution; in addition, the updated trial solutions should not deviate far from the initial trial solution. For DPLR, the good-enough trial solution is exactly the initial trial UC, all we need to do at this time is to make the updated trial solutions close to it. This can be realized by optimizing the initial step size $c_0$ with the following model (which can be easily transformed into a linear program).

$$\begin{cases} \max_{c_0, \hat{\lambda}_1^*} c_0 \\ \text{s.t. } \hat{\lambda}_1^* = \hat{\lambda}^* + c_0 g_{\hat{\lambda}^*}, \\ \left|\beta_{i,t}(\hat{\lambda}_1^*) - \beta_i^0\right| \geq \delta_i, \text{ if } \left|\beta_{i,t}(\hat{\lambda}^*) - \beta_i^0\right| \geq \delta_i, \forall i, t; \\ 0 \leq c_0 \leq u, \hat{\lambda}_{0,t,1}^* \text{ free}, \hat{\lambda}_{l,t,1}^* \geq 0, \forall l, t. \end{cases} \quad (34)$$

where the entries of the constant vector $\delta$ are given scalar thresholds, the constant $g_{\hat{\lambda}^*}$ is given by (28), and the constant $u$ is a given upper bound of the variable $c_0$ representing the initial step size, and the variables $\hat{\lambda}_1^*$ are the updated multipliers at the first iteration expressed as a function of the variable $c_0$ in the first constraint. The second constraint is for keeping the binary solutions $\hat{z}_{i,t}^*(\hat{\lambda}_1^*)$ with bigger $|\beta_{i,t}(\hat{\lambda}_1^*) - \beta_i^0|$ unchanged in the first update of trial UC, and $\delta$ can control the ratio of the unchanged binary solutions to all the binary solutions; if the $|\beta_{i,t}(\hat{\lambda}_1^*) - \beta_i^0|$ of binary solutions $\hat{z}_{i,t}^*(\hat{\lambda}_1^*)$ are smaller than $\delta_i$, $\hat{z}_{i,t}^*(\hat{\lambda}_1^*)$ are allowed to be adjusted. The objective maximizes the initial step size to find the bound $[0, c_0^*]$ of the variable $c_0$ that satisfies the constraints in (34). This problem is feasible since $(0, \lambda^*)$ is feasible. It is found in numerical tests that, with the multipliers calculated with (29) and (34), the method generally requires a few iterations to give the feasible sub-optimal solution $(\bar{z}^*, \bar{p}^*)$, which is of high quality in terms of the generation cost.